\documentclass[12pt,fullpage]{article}

\usepackage{amsfonts}
\usepackage{xr}
\usepackage{graphicx}
\usepackage{amsmath}
\usepackage{amssymb}
\usepackage{epsfig}

\def\R{\mathbb{R}}
\def\N{\mathbb{N}}
\def\C{\mathbb{C}}

\def\epsilon{\varepsilon}
\def\trait (#1) (#2) (#3){\vrule width #1pt height #2pt depth #3pt}
\def\fin{\hfill\trait (0.1) (5) (0) \trait (5) (0.1) (0) \kern-5pt 
\trait (5) (5) (-4.9) \trait (0.1) (5) (0)}

\usepackage{color}

\definecolor{gr}{rgb}   {0.,   0.8,   0. } 
\definecolor{bl}{rgb}   {0.,   0.5,   1. } 
\definecolor{mg}{rgb}   {0.7,  0.,    0.7} 

\newcommand{\Bk}{\color{black}}

\newcommand{\be}{\begin{equation}}
\newcommand{\ee}{\end{equation}}
\newcommand{\baa}{\begin{array}}
\newcommand{\eaa}{\end{array}}
\newcommand{\ba}{\begin{eqnarray}}
\newcommand{\ea}{\end{eqnarray}}

\newtheorem{theo}{\bf Theorem}[section]
\newtheorem{lem}[theo]{\bf Lemma}
\newtheorem{pro}[theo]{\bf Proposition}
\newtheorem{cor}[theo]{\bf Corollary}
\newtheorem{defi}[theo]{\bf Definition}
\newtheorem{rem}[theo]{\bf Remark}

\begin{document}
\author{\begin{tabular}{ccc}
Nadine Badr\footnote{Universit\'e Claude Bernard, Lyon I et CNRS UMR 5208, B\^atiment Doyen Jean Braconnier, 43 boulevard du 11 novembre 1918, 69622 VILLEURBANNE Cedex, FRANCE Email: badr@math.univ-lyon1.fr} &
Emmanuel Russ \footnote{Universit\'e Joseph Fourier, Institut Fourier, UMR 5582, 100 rue des math, BP 74, 38402 Saint-Martin d'H\`eres. E-mail: emmanuel.russ@ujf-grenoble.fr}\\ \\
Universit\'e Claude Bernard & Universit\'e Joseph Fourier 
\end{tabular}}

\title{\Bk A Meyers type regularity result \Bk for approximations of second order elliptic operators by Galerkin schemes}

\maketitle

\noindent{\small{{\bf Abstract.} }} We prove \Bk a Meyers type regularity estimate \Bk for approximate solutions of second order elliptic equations obtained by Galerkin methods. The proofs rely on interpolation results for Sobolev spaces on graphs. Estimates for second order elliptic operators on rather general graphs are also obtained. 

\tableofcontents

\section{Introduction}
Let $\Omega\subset \R^n$ be a bounded $C^2$ domain. Let $A:\Omega\rightarrow {\mathcal M}_n(\R)$ (where ${\mathcal M}_n(\R)$ is the space of $n\times n$ matrices with real coefficients) be a bounded measurable function satisfying the following uniform ellipticity condition: there exists $c>0$ such that, for almost every $x\in \Omega$ and all $\xi\in \R^n$,
\begin{equation} \label{ellipticA}
A(x)\xi\cdot \xi\geq c\left\vert \xi\right\vert^2,
\end{equation}
where, for all $u,v\in \R^n$, $u\cdot v$ is the standard scalar product of $u$ and $v$ and $\left\vert \xi\right\vert$ stands for the Euclidean norm of $\xi$. Given $\Bk f\in W^{-1,2}(\Omega)$, \Bk the Lax-Milgram theorem shows that there exists a unique function $u\in W^{1,2}_0(\Omega)$ such that
\begin{equation} \label{divAgrad}
\mbox{div}(A\nabla u)=f,
\end{equation}
in the sense that, for all function $\varphi\in {\mathcal D}(\Omega)$,
\begin{equation} \label{weak}
\int_{\Omega} A(x)\nabla u(x)\cdot \nabla \varphi(x)dx=-\langle f,\varphi\rangle.
\end{equation}
A classical result of N. Meyers (\cite{meyers}) asserts that there exists $\varepsilon>0$, only depending on $\Omega$ and $A$, such that, for all $p\in (2,2+\varepsilon)$, \Bk if $f\in W^{-1,p}(\Omega)$, \Bk then $u\in W^{1,p}_0(\Omega)$ and 
$$
\Bk \left\Vert u\right\Vert_{W^{1,p}(\Omega)}\leq C_p\left\Vert f\right\Vert_{W^{-1,p}(\Omega)}
$$
\Bk where $C_p>0$ only depends on $\Omega,A$ and $p$. This results still holds when $\Omega$ is Lipschitz (\cite{galmon}). \par
\noindent Let us now focus on the case when $n=2$ and $\Omega\subset \R^2$ is a polygonal domain. Let $A$ as before. We shall be interested in approximate solutions of \eqref{divAgrad} obtained by Galerkin approximations based on spaces of piecewise polynomial functions. Let us briefly present these approximation schemes. By triangulation of $\Omega$, we mean a family ${\mathcal T}$ of triangles such that
$$
\Omega=\bigcup_{T\in {\mathcal T}} T.
$$
Say that ${\mathcal T}$ is admissible if and only if, for all $T,T^{\prime}\in {\mathcal T}$, $T\cap T^{\prime}$ is either empty, or a vertex, or a common side. If ${\mathcal T}$ is an admissible triangulation of $\Omega$, for all $T\in {\mathcal T}$, denote by $h_T$ the diameter of $T$ and by $\rho_T$ the inner diameter of $T$ ({\it i.e.} the diameter of the largest ball included in $T$). Define also
$$
V^1({\mathcal T}):=\left\{v\in C(\Omega);\ v\vert_T\in {\mathcal P}_1\ \forall T\in {\mathcal T}\right\}
$$
where ${\mathcal P}_1$ denotes the space of polynomials with degree less or equal to $1$. Define also $V^1_0({\mathcal T})$ as the subspace of $V^1({\mathcal T})$ made of functions vanishing on $\partial\Omega$. It is plain to see that $V^1({\mathcal T})\subset W^{1,2}(\Omega)$, so that $V^1_0({\mathcal T})\subset W^{1,2}_0(\Omega)$. Another application of the Lax-Milgram theorem gives a unique function $u_{\mathcal T}\in V^1_0({\mathcal T})$ such that
\begin{equation} \label{galerkin}
\int_{\Omega} A(x)\nabla u_{\mathcal T}(x)\cdot \nabla v(x)dx=-\int_{\Omega} f(x)v(x)dx\ \forall v\in V^1_0({\mathcal T})
\end{equation} 
(note that the bilinear form given by the left-hand side of \eqref{galerkin} is still coercive on $V^1_0({\mathcal T})\times V^1_0({\mathcal T})$). Assume now that, for all $h>0$, we have an admissible triangulation ${\mathcal T}_h$ of $\Omega$, such that
$$
\sup_{T\in {\mathcal T}} h_T=h.
$$
Say that the family of triangulations $\left({\mathcal T}_h\right)_{h>0}$ is regular if and only if there exists $\sigma>0$ such that
$$
\frac{h_T}{\rho_T}\leq \sigma\ \forall T\in {\mathcal T}_h,\ \forall h>0.
$$
For all $h>0$, write $u_h$ instead of $u_{{\mathcal T}_h}$. \Bk Then, one always has (see \cite{bs}, Chapter 5):
$$
\lim_{h\rightarrow 0} \left\Vert u_h-u\right\Vert_{W^{1,2}(\Omega)}=0.
$$
In the present work, we establish a Meyers type result for the approximate solutions $u_h$ of \eqref{ellipticA}, as well as an improved convergence result of $u_h$ to $u$ in a Sobolev norm. More precisely, we establish:
\begin{theo} \label{meyersgalerkin}
Assume that $\Omega\subset\R^2$ is a convex polygonal subset and that $A\in L^{\infty}(\Omega)$ is uniformly elliptic, in the sense that $A$ satisfies \eqref{ellipticA}. Then, there exist $\varepsilon>0$ and $C>0$ such that, for all $p\in (2,2+\varepsilon)$ and all $f\in W^{-1,p}(\Omega)$, the solution $u_h\in V^1_0({\mathcal T}_h)$ of \eqref{galerkin} belongs to $W^{1,p}(\Omega)$ and, for all $h>0$,
\begin{equation} \label{estimuh}
\left\Vert u_h\right\Vert_{W^{1,p}(\Omega)}\leq C\left\Vert f\right\Vert_{W^{-1,p}(\Omega)}.
\end{equation}
Moreover, for all $p\in (2,2+\varepsilon)$, 
\begin{equation} \label{convuh}
\lim_{h\rightarrow 0} \left\Vert u_h-u\right\Vert_{W^{1,p}(\Omega)}=0.
\end{equation}
\end{theo}
Theorem \ref{meyersgalerkin} is optimal in the scale of Sobolev spaces in the following sense:
\begin{pro} \label{optimal}
Assume that $\Omega\subset\R^2$ is a convex polygonal subset. Let $\varepsilon>0$. Then, there exists $A\in L^{\infty}(\Omega)$ uniformly elliptic with the following property: for all $p>2+\varepsilon$, there exists $f\in W^{-1,p}(\Omega)$ such that, if $u_h\in V^1_0({\mathcal T}_h)$ is the solution of \eqref{galerkin}, the $(u_h)_{h>0}$ are not uniformly bounded in $W^{1,p}(\Omega)$.
\end{pro}
As an immediate consequence of Theorem \ref{meyersgalerkin} and the usual Sobolev embeddings, we get a version of Theorem \ref{meyersgalerkin} for H\"older spaces:
\begin{cor} \label{holdgalerkin}
Under the assumptions and notations of Theorem \ref{meyersgalerkin}, there exist $\eta>0$ and $C>0$ such that, for all $p\in (2,2+\varepsilon)$ and all $f\in W^{-1,p}(\Omega)$,
\begin{equation} \label{estimholderuh}
\left\Vert u_h\right\Vert_{C^{0,\eta}(\overline{\Omega})}\leq C\left\Vert f\right\Vert_{W^{-1,p}(\Omega)}
\end{equation}
and
\begin{equation} \label{convholderuh}
\lim_{h\rightarrow 0} \left\Vert u_h-u\right\Vert_{C^{0,\eta}(\overline{\Omega})}=0.
\end{equation}
\end{cor}
When $f\in L^2(\Omega)$ and if one furthermore knows that the solution $u$ of \eqref{ellipticA} belongs to $W^{2,2}(\Omega)$ (this is the case for instance if $A\in C^1(\Omega)$), one proves the following conclusion:
\begin{cor} \label{improvedgalerkin}
The assumptions and notations are the same as in Theorem \ref{meyersgalerkin}. Assume that, for all $f\in L^2(\Omega)$, the solution $u$ of \eqref{ellipticA} belongs to $W^{2,2}(\Omega)$. Then, for all $p\in (2,2+\varepsilon)$, there exist $\theta\in (0,1)$ and $C>0$ such that
\begin{equation}Ê\label{convuhsobolev}
\left\Vert u_h-u\right\Vert_{W^{1,p}(\Omega)}\leq Ch^{\theta}\left\Vert f\right\Vert_{L^2(\Omega)}.
\end{equation}
\end{cor}
\Bk
\noindent Let us give a few comments about these results. \par
\noindent \Bk First, at fixed $h$, $u_h$ is piecewise linear. Moreover, since $u_h$ is actually given by a finite system of linear algebraic equations, it is easily seen that there exists $C_h>0$ such that
\begin{equation} \label{rough}
\left\Vert u_h\right\Vert_{W^{1,\infty}(\Omega)}\leq C_h\left\Vert f\right\Vert_{L^2(\Omega)}.
\end{equation}
However, as Proposition \ref{optimal} shows, the constant $C_h$ in \eqref{rough} cannot be bounded when $h\rightarrow 0$ in general. 

\noindent Theorem \ref{meyersgalerkin} and Corollary \ref{holdgalerkin} can be seen as regularity results for approximate solutions of \eqref{divAgrad} obtained by the finite elements method. These results are the counterparts of classical Meyers and De Giorgi type regularity results for the solution of $\mbox{div}(A\nabla u)=f$ (see the results in \cite{galmon,meyers} stated at the very beginning of the paper). Since $\Omega\subset \R^2$, Theorem \ref{meyersgalerkin} yields Corollary \ref{holdgalerkin} at once by means of the classical Sobolev embeddings. \par
\noindent Under the extra assumption that $A$ is symmetric and when the data $f$ belongs tp $L^2(\Omega)$, a version of Corollary \ref{holdgalerkin} was proved in \cite{rey}, Th\'eor\`eme 6.2. Here, we get rid of the symmetry assumption on $A$ and we also derive the convergence result in \eqref{convholderuh}, while only \eqref{estimholderuh} was established in \cite{rey}. To our best knowledge, Theorem \ref{meyersgalerkin} and Proposition \ref{optimal} are completely new. \par

\medskip

\noindent \Bk Using ideas analogous to those involved in the proof of Theorem \ref{meyersgalerkin}, we also establish \Bk Sobolev, \Bk $L^{\infty}$ and $C^{0,\eta}$ estimates for second order uniformly elliptic operators on quite general graphs. A lot of work was devoted to estimates for reversible Markov chains on infinite graphs (see for instance \cite{coulhontrento,cougri,delmo} and the references therein). In these references, the operators under consideration have positive coefficients. In the present work, we deal with more general operators, for which no positivity assumption is made on the coefficients but which are still elliptic (see Section \ref{ellip} below \Bk for the statements and the proofs). \Bk It seems therefore impossible to use techniques such as the Moser iteration or the maximum principle. \par
\medskip

\noindent \Bk Let us describe our strategy for the proof of Theorem \ref{meyersgalerkin} and its corollaries. As in \cite{rey}, the proofs go through the analysis of some second order elliptic operators on graphs. However, we simplify the arguments in a significant way. \Bk Given a triangulation ${\mathcal T}_h$, we define a graph $\Gamma_h$, the vertices (resp. the edges) of which are the vertices (resp. the edges) of ${\mathcal T}_h$, and consider a suitable ``second order elliptic'' operator on $\Gamma_h$, denoted by $L_h$, which approximates $L$ on $\Gamma_h$. We then reduce \eqref{estimuh} to proving that the inverse of $L_h$ is $L^2(\Gamma_h)-W^{1,2+\varepsilon}(\Gamma_h)$ bounded for some $\epsilon>0$.  \par

\medskip

\noindent The proof of this boundedness result for $L_h^{-1}$ is inspired by an analogous perturbation argument for $2m$-th order operators in $\R^{2m}$ (see \cite{asterisque}, Chapter 1). Namely, we first observe that $L_h$ is bounded from $W^{1,p}(\Gamma_h)$ to $W^{-1,p}(\Gamma_h)$ for all $p\in (1,+\infty)$, and that it is an isomorphism from $W^{1,2}(\Gamma_h)$ into $W^{-1,2}(\Gamma_h)$. Then, using the fact that the $W^{1,p}(\Gamma_h)$ and the $W^{-1,p}(\Gamma_h)$ form an interpolation scale for the {\it complex} method (this result is also new and has its own interest), we use a general perturbation result to conclude that $L_h$ is an isomorphism from $W^{1,p}(\Gamma_h)$ onto $W^{-1,p}(\Gamma_h)$ for $p\in (2-\varepsilon,2+\varepsilon)$ for some $\varepsilon>0$, which completes the proof. \par
\noindent We argue similarly to obtain \Bk $W^{1,p}$, \Bk $L^{\infty}$ and $C^{\eta}$ estimates for second order elliptic operators on graphs $\Gamma$, under some geometric assumptions on $\Gamma$. In particular, the volume of a ball with radius $r>0$ has to be controlled from below by $r^2$. \Bk Sobolev estimates are obtained by a perturbation argument relying on interpolation results for Sobolev spaces. Then, the geometric assumptions on $\Gamma$ \Bk imply that $L^2(\Gamma)$ embeds into $W^{-1,p}(\Gamma)$ and $W^{1,p}(\Gamma)$ embeds into $C^{\eta}(\Gamma)$ for $p>2$, which \Bk yields H\"older estimates. \Bk \par

\medskip

\noindent Theorems \ref{meyersgalerkin} and \Bk its corollaries \Bk are stated for open polygonal subsets of $\R^2$. It is plausible that an adaptation of our method would yield an analogue of Theorem \ref{meyersgalerkin} for polyhedral subsets of $\R^3$. However, this would not yield an analogue of Corollary \ref{holdgalerkin}. Obtaining a De Giorgi type regularity result for approximate solutions of \eqref{divAgrad} by a finite elements method is an open challenging problem.\par

\bigskip

\noindent The plan of the paper is as follows. In section \ref{framework}, we present the discrete setting. Section \ref{prel} is devoted to general interpolation results for Sobolev spaces on graphs for the real and the complex method. In Section \ref{sobgaler}, we show Theorem \ref{meyersgalerkin} \Bk and its corollaries, \Bk while estimates for general second order elliptic operators on graphs are stated and established in Section \ref{ellip}. \par
\noindent{\bf Notation: } if $A(f)$ and $B(f)$ are two quantities depending on a function $f$ ranging in a set $E$, say that $A(f)\lesssim B(f)$ if and only if there exists $C>0$ such that, for all $f\in E$,
$$
A(f)\leq CB(f),
$$
and that $A(f)\sim B(f)$ if and only if $A(f)\lesssim B(f)$ and $B(f)\lesssim A(f)$. \par

\medskip

\noindent {\bf Acknowledgements: }we would like to thank Franck Boyer for helpful discussion. 
\section{Presentation of the discrete framework} \label{framework}

Let us give precise definitions of our framework. The following presentation is partly borrowed from \cite{delmo}. In the sequel, $\Gamma$ stands for a non-empty finite or infinite set. When $\Gamma$ is finite, fix a non-empty strict subset of $\Gamma$, called the boundary of $\Gamma$ and denoted by $\partial\Gamma$.
\subsection{The metric on $\Gamma$} 
Let $h_{xy}=h_{yx}\geq 0$ a symmetric weight on $\Gamma\times \Gamma$. 
We assume that $h_{xx}=0$ for all $x\in \Gamma$. If $x,y\in \Gamma$, say that $x\sim y$ if and only 
if $h_{xy}>0$ or $x=y$ .
We define for all $x \in \Gamma$ 
$$h_x:=\sup_{y\sim x}h_{xy}$$
and 
$$h:=\sup_{x,\,y\in \Gamma}h_{xy}.
$$
Denote by $E$ the set of edges in $\Gamma$, {\it i.e.}
\[
E:=\left\{(x,y)\in \Gamma\times \Gamma;\ x\sim y\right\},
\]
and notice that, due to the symmetry of $h$, $(x,y)\in E$ if and only if $(y,x)\in E$. In the sequel, 
we always assume that $\Gamma$ is locally uniformly finite, which means 
that there exists $N\in \N^{\ast}$ such that
\begin{equation}\label{uni}
\sup_{x\in \Gamma} \sharp \{y\in \Gamma; \, y\sim x\} \leq N
\end{equation}
(here and after, $\sharp A$ denotes the cardinal of any subset $A$ of $\Gamma$).\par
\noindent We also assume that $h$ is  locally uniformly controlled with a constant $C_W\geq 1$, which means that for all $x,\,y, \, z\in \Gamma$ such that $x\sim y$  with $x\neq y$ and   $x\sim z$  with $x\neq z$, we have
\begin{equation} \label{controlh}
\frac{1}{C_W}\leq \frac{h_{xy}}{h_{xz}}\leq C_{W}.
\end{equation}
\Bk For $x,y\in \Gamma$, a path $C$ joining $x$ to $y$ is a finite sequence of vertices $x_0=x,...,x_N=y$ such that, for all $0\leq i\leq N-1$, 
$x_i\sim x_{i+1}$. 
By definition, the length of such a path is
$$l(C):=\sum_{k=0}^{N-1}h_{x_k x_{k+1}}.$$
Assume that $\Gamma$ is connected, which means that, for all $x,y\in \Gamma$, there exists a path joining 
$x$ to $y$. For all $x,y\in \Gamma$, the distance between $x$ and $y$, 
denoted by $d(x,y)$, is defined as the infimum of the lengths of all paths joining $x$ and $y$. For all $x\in \Gamma$ and all $r>0$, let 
$B(x,r):=\left\{y\in \Gamma,\ d(y,x)< r\right\}$. 
\begin{rem} \label{singleton}
If $x\in \Gamma$ and $r>0$ are such that $r\leq \frac{h_x}{C_W}$, then $B(x,r)=\{x\}$. Indeed, let $y\in B(x,r)$, so that $d(y,x)<r\leq \frac{h_x}{C_W}$. If $y\neq x$, there exists $z\sim x$ such that 
$$
h_{xz}\leq d(y,x)<\frac{h_x}{C_W}\leq h_{xz},
$$
which is impossible.
\end{rem}

If $B=B(x,r)$ is a ball, set $\alpha B:=B(x,\alpha r)$ for all $\alpha>0$.
\subsection{A measure on the graph and the edges}
\noindent We also consider measures on $E$ and $\Gamma$. Assume that, for all $x\sim y\in \Gamma$, we are given a coefficient $\mu_{xy}> 0$ such that $\mu_{xy}=\mu_{yx}$. Assume also that there exists $C>0$ such that, for all $x,y,z\in \Gamma$ satisfying $x\sim y$ and $x\sim z$,
\begin{equation} \label{comparmu}
\frac{\mu_{xy}}{\mu_{xz}}\leq C.
\end{equation}
For all $x\in \Gamma$, set
\begin{equation} \label{defm}
m(x):=\sum_{y\sim x} \mu_{xy}.
\end{equation}
If $A\subset \Gamma$, let
$$
m(A):=\sum_{x\in A} m(x).
$$

\noindent For all $1\leq p<+\infty$, say that a function $f$ on $\Gamma$ belongs to $L^p(\Gamma,m)$ (or $L^p(\Gamma)$) if and only if
\[
\left\Vert f\right\Vert_p:=\left(\sum\limits_{x\in \Gamma}\left\vert f(x)\right\vert^pm(x)\right)^{1/p}<+\infty.
\]
Say that $f\in L^{\infty}(\Gamma,m)$ (or $L^{\infty}(\Gamma)$) if and only if
\[
\left\Vert f\right\Vert_{\infty}:=\sup\limits_{x\in \Gamma} \left\vert f(x)\right\vert<+\infty.
\]
When $\Gamma$ is finite and $1\leq p\leq +\infty$, say that $f\in L^p_0(\Gamma)$ if and only if  $f\in L^p(\Gamma)$ and $f(x)=0$ for all $x\in \partial\Gamma$.
\begin{rem} \label{denselp}
It is an elementary observation that the space of functions on $\Gamma$ with finite support is dense in $L^p(\Gamma,m)$ for all $1\leq p<+\infty$.
\end{rem}

\medskip

\noindent For all $1\leq p<+\infty$, say that a function $F$ on $E$ belongs to $L^p(E,\mu)$ (or $L^p(E)$) if and only if $F$ is antisymmetric (which means that $F(x,y)=-F(y,x)$ for all $(x,y)\in E$) and
\[
\left\Vert F\right\Vert_p:=\left(\sum\limits_{(x,y)\in E}\left\vert F(x,y)\right\vert^p\mu_{xy}\right)^{1/p}<+\infty.
\]
Say that $F\in L^{\infty}(E,\mu)$ (or $L^{\infty}(E)$) if and only $F$ is antisymmetric and
\[
\left\Vert F\right\Vert_{\infty}:=\sup\limits_{(x,y)\in E} \left\vert F(x,y)\right\vert<+\infty. 
\]
\subsection{Differential and gradient}
In the spirit of \cite{badrruss}, for all function $f:\Gamma\rightarrow \R$ and all $(x,y)\in E$, define
$$
df(x,y)=\left\{
\begin{array}{ll}
\displaystyle \frac{f(y)-f(x)}{h_{xy}} & \mbox{ if }x\neq y,\\
0 & \mbox{ if }x=y.
\end{array}
\right.
$$
Say that $d$ is the differential operator on $\Gamma$. \par \Bk

\medskip

\noindent We also define the sublinear operator ``length of the gradient''  by
\[
\nabla f(x)=\frac {1}{h_{x}}\left(\sum\limits_{y\sim x} \left\vert 
f(y)-f(x)\right\vert^{2}\right)^{1/2}
\] 
for all function $f$ on $\Gamma$ and all $x\in \Gamma$. It is a consequence of \eqref{uni}, \eqref{controlh}, \eqref{comparmu} and \eqref{defm} that, for all function $f:\Gamma\rightarrow \R$ and all $p\in [1,+\infty]$, 
$$
\left\Vert df\right\Vert_{L^p(E,\mu)}\sim \left\Vert \nabla f\right\Vert_{L^p(\Gamma,m)}.
$$

\subsection{Sobolev and H\"older spaces}
Let $1\leq p\leq +\infty$. Say that a scalar-valued function $f$ on $\Gamma$ belongs to the (inhomogeneous) Sobolev space $W^{1,p}(\Gamma)$ (see also \cite{badrruss,gt,ostr}) if and only if
\[
\left\Vert f\right\Vert_{W^{1,p}(\Gamma)}:=\left\Vert f\right\Vert_{L^p(\Gamma)}+\left\Vert \nabla f\right\Vert_{L^p(\Gamma)} <+\infty.
\]
In the case when $\Gamma$ is finite, denote by $W^{1,p}_0(\Gamma)$ the subspace of $W^{1,p}(\Gamma)$ made of functions $f$ such that $f=0$ on $\partial \Gamma$.\par
\noindent It is then routine to check that both $W^{1,p}(\Gamma)$ and $W^{1,p}_0(\Gamma)$ are Banach spaces. \par

\medskip

\noindent Finally, if $1<p<+\infty$, as in the Euclidean case, we  define $W^{-1,p}(\Gamma)$ as the dual space of $W^{1,p^{\prime}}(\Gamma)$ (when $\Gamma$ is infinite) or of $W^{1,p^{\prime}}_0(\Gamma)$ (when $\Gamma$ is finite), equipped with its natural norm. Here and after, $p^{\prime}$ is defined by $\frac 1p+\frac 1{p^{\prime}}=1$.
\begin{rem} \label{densesobol}
It is easy to check that, for all $1\leq p<+\infty$, the space of functions with finite support is dense in $W^{1,p}(\Gamma)$.
\end{rem}
\bigskip

\noindent For all $\mu>0$, we define the discrete homogeneous (resp. inhomogeneous) H\"older spaces of exponent $\mu$ as 
$$ \dot{C}^{\mu}(\Gamma):=\{ u\in \mathbb{R}^{\Gamma};\, |u|_{\dot{C}^{\mu}(\Gamma)}= \sup_{x\neq y}\frac{|u(x)-u(y)|}{d(x,y)^{\mu}}< +\infty \}
$$
and
$$ C^{\mu}(\Gamma):=\{ u\in \mathbb{R}^{\Gamma};\, \|u\|_{C^{\mu}(\Gamma)}= \|u\|_{L^{\infty}(\Gamma)}+|u|_{\dot{C}^{\mu}(\Gamma)}<+\infty \}. 
$$

\subsection{Geometric assumptions on $\Gamma$}
Let us now introduce some geometric assumptions on $(\Gamma,m)$.

\begin{itemize}
\item[1.] Say that $(\Gamma,m)$ satisfies the doubling property if there exists $C>0$ such that, for all $x\in \Gamma$ and all $r>0$,  
\begin{equation} \label{D} \tag{$D$}
V(x,2r)\leq CV(x,r).
\end{equation}
Note that this assumption implies that there exist $C,D>0$ such that, for all $x\in \Gamma$, all $r>0$ and all $\theta>1$,
\begin{equation} \label{Dbis}
V(x,\theta r)\leq C\theta^DV(x,r).
\end{equation}
\begin{rem} \label{infinite}
Observe also that, when $\Gamma$ is infinite, it is also unbounded (since it is locally uniformly finite) so that, if (\ref{D}) holds, then $m(\Gamma)=+\infty$ (see \cite{martellthesis}). 
\end{rem}

\item[2.] Say that $(\Gamma,m)$ satisfies $(L_{\sigma})$ if there exist $\sigma>0,\, c>0$ such that, for all $x\in \Gamma$ and all $r>0$,  
\begin{equation} \label{L} \tag{$L_{\sigma}$}
V(x,r)\geq cr^{\sigma}.
\end{equation}

\item[3.] If $1\leq q<+\infty$,  say that $(\Gamma,m)$ satisfies a scaled $L^{q}$ Poincar\'e 
inequality on balls (this inequality will be denoted by $(P_{q})$ in 
the sequel) if there exists 
$C>0$ such that, for any $x\in \Gamma$, any $r>0$ and any function 
$f$ on $\Gamma$,
\begin{equation} \label{poincarel2} \tag{$P_{q}$}
\sum_{y\in B(x,r)} \left\vert f(y)-f_B\right\vert^{q}m(y)\leq 
Cr^{q}\sum_{y\in B(x,r)} \left\vert \nabla f(y)\right\vert^{q}m(y),
\end{equation}
where
\[
f_B=\frac 1{V(B)} \sum_{x\in B} f(x)m(x)
\]
is the mean value of $f$ on $B$. 
\end{itemize}
When $\Gamma$ is finite, we will consider local versions of these properties. Let $r_0>0$.
Say that $\Gamma$ satisfies $(D_{loc})$ if there exists $C>0$ such that, for all $0<r<r_0$ and all $x\in \Gamma$,
$$
V(x,2r)\leq CV(x,r).
$$
Say that $\Gamma$ satisfies $(L_{\sigma,loc})$ if there exists $C>0$ such that, for all $x\in \Gamma$ and $0<r\leq  r_0$,
$$
V(x,r)\geq Cr^{\sigma}.
$$
Let $1\leq q<+\infty$. Say that $\Gamma$ satisfies $(P_{q,loc})$ if there exists $C>0$ such that, for all $x\in \Gamma$, all $0<r<r_0$ and any function $f$ on $\Gamma$,
$$
\sum_{y\in B(x,r)} \left\vert f(y)-f_B\right\vert^{q}m(y)\leq 
Cr^{q}\sum_{y\in B(x,r)} \left\vert \nabla f(y)\right\vert^{q}m(y).
$$

\section{Interpolation and embeddings of Sobolev spaces}\label{prel}
\subsection{Real and complex interpolation of Sobolev spaces: the case of unbounded graphs}\label{intrecounbounded}
Let us first focus on the case when $\Gamma$ is infinite, and therefore unbounded. We will make constant use of the following notation here and in the sequel:
\begin{defi} \label{q0}
When $(P_q)$ holds for some $1\leq q<\infty$, we define  $q_{0}=\inf\left\lbrace q \in [1,\infty):(P_{q})\textrm{ holds}\right\rbrace$.
\end{defi}
Note that if $(P_q)$ holds for some $q>1$, then $q_0<q$ (see \cite{kz}), and that $(P_r)$ holds for all $q_0<r<+\infty$.\par

\bigskip

\noindent Let us first state the following real interpolation theorem for Sobolev spaces:
\begin{theo} \label{interpolation}
Let $q\in [1,+\infty)$ and assume that (\ref{D}) and $(P_q)$ hold. Then, for all $1\leq r\leq q<p<+\infty$, $W^{1,p}(\Gamma)=\left(W^{1,r}(\Gamma),W^{1,\infty}(\Gamma)\right)_{1-\frac rp,p}$.
\end{theo}
The corresponding result for Sobolev spaces on Riemannian manifolds was proved in \cite{badr1}, Theorem 1.1. We established in \cite{badrruss}, Theorem 1.18, a statement very similar to Theorem \ref{interpolation} for homogeneous Sobolev spaces, and the proof, which we skip here, is analogous to the one of Theorem 1.18 in \cite{badrruss}, even if the definition of the gradient is slightly different. The only important assumptions to be able to interpolate are \eqref{D} and $(P_q)$, see Theorem 7.11 in \cite{badr1}. \par
\medskip

\noindent As an immediate corollary, we obtain:
\begin{cor}[The reiteration theorem]\label{RH} Assume that $\Gamma$ satisfies (\ref{D}), $(P_{q})$ for some $1\leq q<+\infty$. Let $q_{0}<p_{1}<p<p_{2}\leq +\infty$, with $\displaystyle \frac 1p=\frac{1-\theta}{p_1}+\frac{\theta}{p_2}$, then $$
W^{1,p}(\Gamma)=\left(W^{1,p_1}(\Gamma),W^{1,p_2}(\Gamma)\right)_{\theta,p}.
$$
If $q_0=1$, then one can also take $1=p_1<p<p_2\leq +\infty$.
\end{cor}
\noindent Let us now turn to the complex interpolation of Sobolev spaces. We first recall the following complex reiteration theorem:
\begin{theo} \label{cwikel}\cite{bl,cwikel}
For any compatible couple of Banach spaces $(A_1,A_2)$ we have
$$
[(A_1,A_2)_{\lambda_1,p_2}, (A_1,A_2)_{\lambda_2,p_2}]_{\alpha}=(A_1,A_2)_{\beta,p}
$$
for all $\lambda_1,\lambda_2$ and  $\alpha$ in $(0,1)$ and all $p_1, p_2$ in $[1,+\infty]$, except for the case $p_1=p_2=\infty$.
Here $\beta$ and $p$ are given by $\beta=(1-\alpha)\lambda_1+\alpha \lambda_2$ and $\frac{1}{p}=\frac{1-\alpha}{p_1}+\frac{\alpha}{p_2}.$
\end{theo}
From this, we deduce a complex interpolation result for  Sobolev spaces:
\begin{cor} \label{complex1} Assume that $\Gamma$ satisfies (\ref{D}) and $(P_{q})$ for some $1\leq q<+\infty$. Let $q_0<p_1<p<p_2<+ \infty$ where $q_0$ is given by Definition \ref{q0}. Then for $\alpha= \frac{\frac{1}{p_1}-\frac{1}{p}}{ \frac{1}{p_1}-\frac{1}{p_2}}= \frac{p_1(p-p_1)}{p(p_2-p_1)}$ one has
$$\left[W^{1,p_1}(\Gamma), W^{1,p_2}(\Gamma)\right]_{\alpha}=W^{1,p}(\Gamma).$$
\end{cor}
\textbf{Proof:}  We apply Theorem \ref{cwikel} with $A_1=W^{1,1}(\Gamma), \, A_2=W^{1,\infty}(\Gamma)$, $\lambda_1=1-\frac{1}{p_1}$, $\lambda_2=1-\frac{1}{p_2}$  and $\beta= 1-\frac{1}{p}$ .
We obtain 
$$\left[W^{1,p_1}(\Gamma), W^{1,p_2}(\Gamma)\right]_{\alpha}= (W^{1,1}(\Gamma),W^{1,\infty}(\Gamma))_{1-\frac{1}{p},p}
$$
and the result follows from Theorem \ref{interpolation}  which yields
$$(W^{1,1}(\Gamma),W^{1,\infty}(\Gamma))_{1-\frac{1}{p},p}= W^{1,p}(\Gamma)$$
since $p>q_0$.
\hfill\fin\par

Now for dual spaces, we recall the following theorem:
\begin{theo}\label{dual} (\cite{bl}, Theorem 3.7.1)
Let $A_1, A_2$ a compatible couple of  Banach spaces. Then for all $\theta \in (0,1)$ and  $1\leq s \leq \infty$, we have
$$
[(A_1,A_2)_{\theta,s}]^*=(A_1^*,A_2^*)_{\theta,s'}.$$
\end{theo}

Combining this Theorem with Theorem \ref{interpolation}, we deduce:

\begin{theo}\label{interpolationdu} Assume that $\Gamma$ satisfies (\ref{D}) and $(P_{q})$ for some $1\leq q<+\infty$. Then, for all $1\leq r_1<r<r_2\leq \infty$ with $r^{\prime}>q_0$,  we have
$$
\left(W^{-1,r_1}(\Gamma),W^{-1,r_2}(\Gamma)\right)_{\theta,r}=W^{-1,r}(\Gamma)
$$
where $\theta \in (0,1)$ such that $\frac{1}{r}=\frac{1-\theta}{r_1}+\frac{\theta}{r_2}$.
\end{theo}

\textbf{Proof:} 
We have $1\leq r'_2<r'<r'_1\leq \infty$ with $r^{\prime}>q_0$ and $\frac{1}{r'}=\frac{\theta}{r'_2}+\frac{1-\theta}{r'_1}.$ Then

$$
\begin{array}{l}
\displaystyle [(W^{1,r'_2}(\Gamma),W^{1,r'_1}(\Gamma))_{1-\theta,r'}]^*=(W^{1,r'}(\Gamma))^*\\
\displaystyle =W^{-1,r}(\Gamma)= (W^{-1,r_2}(\Gamma),W^{-1,r_1}(\Gamma))_{1-\theta,r}=(W^{-1,r_1}(\Gamma),W^{-1,r_2}(\Gamma))_{\theta,r}.
\end{array}
$$
\hfill\fin\par
\begin{theo} \label{complex2} Assume that $\Gamma$ satisfies (\ref{D}), $(P_{q})$ for some $1\leq q<+\infty$.
Let  $1<p_1<p<p_2<q'_0$, where $q_0$ is given by Definition \ref{q0}.  Then  
$$
\left[W^{-1,p_1}(\Gamma),W^{-1,p_2}(\Gamma)\right]_{\alpha}= W^{-1,p}(\Gamma)$$
where $\beta=1-\frac{1}{p}$ and $\alpha= \frac{\frac{1}{p_1}-\frac{1}{p}}{ \frac{1}{p_1}-\frac{1}{p_2}}$.
\end{theo}

\noindent \textbf{Proof:}
We have $p_1<p<p_2<q'_0$. Then  Theorem \ref{interpolationdu} shows that
$$
(W^{-1,1}(\Gamma),W^{-1,\infty}(\Gamma))_{1-\frac{1}{p_1},p_1}= W^{-1,p_1}(\Gamma),\  (W^{-1,1}(\Gamma),W^{-1,\infty}(\Gamma))_{1-\frac{1}{p_2},p_2}
= W^{-1,p_2}(\Gamma)
$$
and
$$
(W^{-1,1}(\Gamma),W^{-1,\infty}(\Gamma))_{1-\frac{1}{p},p}= W^{-1,p}(\Gamma).
$$
This, with Theorem \ref{cwikel}, yields
$$ \left[W^{-1,p_1}(\Gamma), W^{-1,p_2}(\Gamma)\right]_{\alpha}=(W^{-1,1}(\Gamma),W^{-1,\infty}(\Gamma))_{\beta,p}=W^{-1,p}(\Gamma).$$
\hfill\fin
\subsection{Real and complex interpolation of Sobolev spaces: the case of bounded graphs} \label{intrecobounded}
Assume now that the graph $\Gamma$ is bounded.
\begin{theo} \label{interpolationb}
Let $q\in [1,+\infty)$ and assume that  $(D_{loc})$ and $(P_{q,loc})$ hold. Then, for all $1\leq r\leq q<p<+\infty$, $W^{1,p}_0(\Gamma)=\left(W^{1,r}_0(\Gamma),W^{1,\infty}_0(\Gamma)\right)_{1-\frac rp,p}$.
\end{theo}
\textbf{Proof}
We refer to \cite{badr1} and \cite{badrruss} for the proof. We just mention the differences between the proof in \cite{badr1,badrruss} and the present situation. We have to estimate the functional $K$ of real interpolation. Let $$\Omega=\{x\in \Gamma; \mathcal{M}( (|\nabla f| +|f|)^q)(x)>\alpha^q(t)\}$$
where $\mathcal{M}$ is the  uncentered Hardy-Littlewood maximal function, that is
$$
{\mathcal M}h(x):=\sup_{B\ni x} \frac 1{V(B)}\sum_{y\in B} \left\vert h(y)\right\vert m(y)
$$
and $\alpha(t) = \mathcal{M}( (|\nabla f| +|f|)^q)^{*\frac{1}{q}}(t)$.We recall that
 for a function $g$ on $\Gamma$, $g^*$ denotes its decreasing rearrangement function. \par
\noindent If $\Omega =\Gamma$ (note that this may happen since $m(\Gamma)<+\infty$), this is the easy case, we argue as in \cite{badr1}, Section 3.2.2. \par
\noindent Otherwise, we write $\Gamma$ as the  union of balls $B^{i}$, of radius $\rho$ small enough (namely $0<\rho<r_0$) and  having the bounded overlap property. We take $(\varphi^{i})_i$ a partition of unity subordinated to the covering $(B^{i})_i$ and write $f=\sum_if\varphi^{i}=\sum_i f^{i}$. We  estimate  the functional $K$ as in \cite{badr1}, Section 3.2.2, using the Calder\'{o}n-Zygmund decomposition applied to $f^{i}$. 
Note that here  the functions $\left(b^{i}_j\right)_j$ of the decomposition should belong to $W^{1,r}_0(\Gamma)$, that is $b^{i}_j=0$ on $\partial \Gamma$. For this, we take $b_{j}^i= f^{i}\chi^{i}_j$ if $B^{i}_j\cap \partial \Gamma \neq \emptyset$ and, if not, we define as usual $b^{i}_j=(f^i-(f^i)_{B^{i}_j})\chi^i_j$.  We characterize then $K$ as in \cite{badr1}. Integrating $K$, we get the interpolation result.
\hfill\fin\par

\bigskip

As in section \ref{intrecounbounded}, from this theorem, we deduce:
\begin{cor} \label{complexb} Assume that $\Gamma$ satisfies $(D_{loc})$ and $(P_{q,loc})$ for some $1\leq q<+\infty$. Let $q_0<p_1<p<p_2<+ \infty$ where $q_0$ is defined by \eqref{q0}. Then for $\alpha= \frac{\frac{1}{p_1}-\frac{1}{p}}{ \frac{1}{p_1}-\frac{1}{p_2}}= \frac{p_1(p-p_1)}{p(p_2-p_1)}$ one has
$$\left[W^{1,p_1}_0(\Gamma), W^{1,p_2}_0(\Gamma)\right]_{\alpha}=W^{1,p}_0(\Gamma).$$
\end{cor}

As far as dual spaces are concerned, we have:
\begin{theo}\label{interpolationdub} Assume that $\Gamma$ satisfies $(D_{loc})$ and $(P_{q,loc})$ for some $1\leq q<+\infty$. Then, for all $1\leq r_1<r<r_2\leq \infty$ with $r^{\prime}>q_0$,  we have
$$
\left(W^{-1,r_1}(\Gamma),W^{-1,r_2}(\Gamma)\right)_{\theta,r}=W^{-1,r}(\Gamma)
$$
where $\theta \in (0,1)$ such that $\frac{1}{r}=\frac{1-\theta}{r_1}+\frac{\theta}{r_2}$.
\end{theo}
Finally,
\begin{theo} \label{complexb2} Assume that $\Gamma$ satisfies $(D_{loc})$, $(P_{q,loc})$ for some $1\leq q<+\infty$.
Let  $1<p_1<p<p_2<q'_0$.  Then  
$$
\left[W^{-1,p_1}(\Gamma),W^{-1,p_2}(\Gamma)\right]_{\alpha}= W^{-1,p}(\Gamma)$$
where $\beta=1-\frac{1}{p}$ and $\alpha= \frac{\frac{1}{p_1}-\frac{1}{p}}{ \frac{1}{p_1}-\frac{1}{p_2}}$
\end{theo}

\subsection{Sobolev embeddings} \label{sobemb}
Assuming (\ref{D}), $(P_{q})$ and $(L_{\sigma})$ for some $\sigma>0$, we obtain Sobolev embeddings analogous to the Euclidean situation. Namely:
\begin{pro} \label{sobolev} 
Assume that $\Gamma$ satisfies (\ref{D}), $(P_{q})$ and $(L_{\sigma})$ for some $\sigma>0$.
Then 
$$ \|f\|_{L^{p^*}(\Gamma)}\lesssim \|\nabla f\|_{L^p(\Gamma)} $$
for all $q_0<p< \sigma$ where $p^*=\frac{\sigma p}{\sigma-p}$.
\end{pro}
This result is proved in \cite{saloffcoste}, Corollary 3.3.3, p.75.
\begin{pro} \label{holder} Assume that $\Gamma$ satisfies (\ref{D}), $(P_{q})$ and $(L_{\sigma})$ for some $\sigma>0$. Then, for all $p>\max(q_0,\sigma)$,
$$W^{1,p}(\Gamma)\hookrightarrow C^{\eta}(\Gamma)$$
with $\eta:= 1-\frac{\sigma}{p}$, where the symbol $\hookrightarrow$ means that the embedding is continuous.
\end{pro}

\textbf{Proof:} Let $u \in W^{1,p}(\Gamma)$, $x,y\in \Gamma$ with $x\neq y$ and set $r:=d(x,y)$. Define $B:=B(x,r)$, so that $B\subset B(x,2r)\cap B(y,2r)$. One has
\begin{equation} \label{holderu}
\left\vert u(x)-u(y)\right\vert\leq \left\vert u(x)-u_B\right\vert+ \left\vert u(y)-u_B\right\vert,
\end{equation}
and we estimate each term of the right-hand side of \eqref{holderu}. \par
\noindent Let us first focus on the second term. Using Remark \ref{singleton}, write
\begin{align*}
|u(y)-u_{B}| &\leq \sum_{k=-\infty}^{0}|u_{B(y,2^{k}r)}-u_{B(y,2^{k+1}r)}|+|u_{B(y,2r)}-u_{B}|
\\
&\leq \sum_{k=-\infty}^{0}\frac{1}{V(y,2^kr)}\sum_{z\in B(y,2^kr)}|u(z)-u_{B(y,2^{k+1}r)}|m(z)+\frac{1}{V(B)}\sum_{z\in B}|u(z)-u_{B(y,2r)}|m(z)
\\
&\leq \sum_{k=-\infty}^{0}\frac{V(y,2^{k+1}r)}{V(y,2^{k}r)}\frac{1}{V(y,2^{k+1}r)}\sum_{z\in B(y,2^{k+1}r)}|u(z)-u_{B(y,2^{k+1}r)}|m(z)
\\
&+\frac{C}{V(y,2r)}\sum_{z\in B(y,2r)}|u(z)-u_{B(y,2r)}|m(z)\\
&\leq C\sum_{k=-\infty}^{0}2^kr\left(\frac{1}{V(y,2^{k+1}r)}\sum_{z\in B(y,2^{k+1}r)}|\nabla u(z)| ^p m(z)\right)^{\frac{1}{p}}+Cr \left(\frac 1{V(y,2r)} \sum_{z\in B(y,2r)} \left\vert \nabla u(z)\right\vert^p\right)^{\frac 1p} 
\\
&\leq C\left(\sum_{k=-\infty}^{0}\left(2^{k}r\right)^{1-\frac{\sigma}{p}}\right)\| \,\nabla u\,\|_{L^p(\Gamma)} +Cr^{1-\frac{\sigma}p} \left\Vert \nabla u\right\Vert_{L^p(\Gamma)} 
\\
&\leq  C r^{1-\frac{\sigma}{p}} \|\,\nabla u\,\|_{L^p(\Gamma)} =C d(x,y)^{1-\frac{\sigma}{p}} \left\Vert \nabla u\right\Vert_{L^p(\Gamma)},
\end{align*}
where we used $(D)$ and $(P_p)$ in the third and fourth inequalities and $(L_{\sigma})$ in the fifth one. \par
\noindent Arguing similarly, one obtains
\begin{equation} \label{holderux}
|u(x)-u_{B}| \leq C d(x,y)^{1-\frac{\sigma}{p}} \left\Vert \nabla u\right\Vert_{L^p(\Gamma)}.
\end{equation}
Thus $u \in \dot{C}^{\eta}(\Gamma)$ and $|u|_{\dot{C}^{\eta}(\Gamma)}\leq C \| \,\nabla u\,\|_{L^p(\Gamma)}$.\par
\noindent It remains to prove that $u\in L^{\infty}(\Gamma)$. For all $x\in \Gamma$,
write
\begin{equation}
\label{IT}
|u(x)|\leq |u(x)-u_{B(x,1)}|+|u_{B(x,1)}-u(y)|+ |u(y)|
\end{equation}
for all $y\in B(x,1)$.\par
\noindent Applying \eqref{holderux} with $r=1$, one obtains
 \begin{equation}\label{Eu}
 |u(x)-u_{B(x,1)}|\leq C\left\Vert \nabla u\right\Vert_{L^p(B(x,1))}.
 \end{equation}
Then, taking the average over $B(x,1)$ in (\ref{IT}) yields 
\begin{align*}
|u(x)| &\leq |u(x)-u_{B(x,1)}|+\frac{1}{V(x,1)} \sum_{y\in B(x,1)}|u(y)-u_{B(x,1)}|m(y)+\frac{1}{V(x,1)} \sum_{y\in B(x,1)}|u(y)|m(y) 
\\
&\leq C\left(1+\frac 1{V(x,1)^{\frac{1}{p}}}\right)\left\Vert \,\nabla u\,\right\Vert_{L^p(B(x,1))}+\frac{C}{V(x,1)^{\frac{1}{p}}} \| u\|_{L^p(B(x,1))}
\\
&\leq C \|  u\|_{W^{1,p}(\Gamma)}.
\end{align*}
We used  (\ref{Eu}), $(P_p)$ and $(L_{\sigma})$ in the last inequality. Thus $\|u\|_{L^{\infty}(\Gamma)}\leq C  \|  u\|_{W^{1,p}(\Gamma)}$ and therefore $u\in C^{\eta}(\Gamma)$  with  $\|  u\|_{ C^{\eta}(\Gamma)}\leq C \|  u\|_{W^{1,p}(\Gamma)}$.
\hfill\fin\par

\section{\Bk Meyers type regularity for Galerkin schemes \Bk} \label{sobgaler}
This section is devoted to the proof of Theorem \ref{meyersgalerkin}, Proposition \ref{optimal}, Corollary \ref{holdgalerkin} and Corollary \ref{improvedgalerkin}. 
\subsection{Definition of a graph and an operator}
Let $\Omega\subset \R^2$ be a convex polygonal domain and $\left({\mathcal T}_h\right)_{h>0}$ be a regular family of triangulations of $\Omega$ as in the statement of \Bk Theorem \ref{meyersgalerkin}. \Bk
\noindent Let $\Gamma_h$ be the set of vertices of ${\mathcal T}_h$. Define $\partial\Gamma_h:=\Gamma_h\cap \partial\Omega$. For all $x,y\in \Gamma_h$, say that $x\sim y$ if and only if $x=y$ or there exists a triangle $T$ in ${\mathcal T}_h$ such that the segment joining $x$ and $y$ is an edge of $T$. For all $x\sim y\in \Gamma_h$, set 
$$
h_{xy}=\left\vert x-y\right\vert,
$$
where we recall that $\left\vert u\right\vert$ denotes the Euclidean norm of $u\in \R^2$. Note that, since the triangulation is admissible, \eqref{uni} and \eqref{controlh} hold. For all $x\sim y$, define also
$$
\mu_{xy}:=h_{xy}^2,
$$
so that \eqref{comparmu} holds and
$$
m(x)\sim h_x^2
$$
for all $x\in \Gamma_h$. 
\begin{lem} \label{geomgraph}
The graph $\Gamma_h$, endowed with $h_{xy}$ and $\mu_{xy}$ satisfies $(D_{loc})$, $(P_{2,loc})$ and $(L_{2,loc})$, with constants \Bk only depeding on $\Omega$ and $\sigma$ and independent of $h$. \Bk
\end{lem}
{\bf Proof: } The proof uses essentially the regularity of the graph $\Gamma_h$ and Poincar\'e inequality $(P_{2,loc})$ on $\Omega$. A detailed proof can be found in \cite{rey}, Chapter 7, section 7.2.
 \hfill\fin\par
\medskip

\Bk \noindent For all $x\in \Gamma_h\cap\Omega$, define $\varphi_x$ as the unique function in $V^1_0({\mathcal T}_h)$ such that 
$$
\varphi_x(y)=\delta_{xy}
$$
for all $y\in \Gamma_h$. The functions $(\varphi_x)_{x\in \Gamma_h\cap\Omega}$ form a basis of $V^1_0({\mathcal T}_h)$. \par
\noindent Let $R_h$ be the ``reconstruction'' operator defined by
$$
R_h\widetilde{u_h}:=\sum_{x\in \Gamma_h\cap\Omega} \widetilde{u_h}(x)\varphi_x
$$
for all function $\widetilde{u_h}$ on $\Gamma_h$ vanishing on $\partial\Gamma_h$. In other words, $R_h\widetilde{u_h}$ is the function in $V^1_0({\mathcal T}_h)$ which coincides with $\widetilde{u_h}$ at the vertices of ${\mathcal T}_h$. We claim that the norms of $\widetilde{u_h}$ and $R_h\widetilde{u_h}$ in the various functional spaces introduced before are equivalent:
\begin{pro} \label{equivnorm}
Let $p\in (1,+\infty)$ and $\eta\in (0,1)$. For all function $\widetilde{u_h}$ on $\Gamma_h$ vanishing on $\partial\Gamma_h$,
\begin{itemize}
\item[$1.$]
$\left\Vert \widetilde{u_h}\right\Vert_{L^p(\Gamma_h)}\sim \left\Vert R_h\widetilde{u_h}\right\Vert_{L^p(\Omega)}$,
\item[$2.$]
$\left\Vert \nabla \widetilde{u_h}\right\Vert_{L^{p}(\Gamma_h)}\sim \left\Vert \nabla R_h\widetilde{u_h}\right\Vert_{L^p(\Omega)}$,
\item[$3.$]
$\left\Vert \widetilde{u_h}\right\Vert_{C^{\eta}(\Gamma_h)}\sim \left\Vert R_h\widetilde{u_h}\right\Vert_{C^{0,\eta}(\overline{\Omega})}$.
\end{itemize}
The implicit constants in these equivalences of norms \Bk only depend on $\sigma$ and are independent of $h$. \Bk
\end{pro}
{\bf Proof: } It is a straightforward consequence of the definitions of the norms, of $R_h$, the properties of the triangulation: namely the regularity of $\Gamma_h$ and the fact that it is locally unifomly bounded, and finally the Poincar\'e inequality on a bounded open convex subset of $\R^2$.  A detailed proof can be found in Appendix E of \cite{rey}.\hfill\fin\par
\Bk \noindent Let us now define a maximal accretive operator on $L^2_0(\Gamma_h)$ associated to the Galerkin scheme. For all $\widetilde{u_h},\widetilde{v_h}\in W^{1,2}_0(\Gamma_h)$, define
$$
Q_h\left(\widetilde{u_h},\widetilde{v_h}\right):=\int_{\Omega} A(x)\nabla\left(R_h\widetilde{u_h}\right)(x)\cdot \nabla\left(R_h\widetilde{v_h}\right)(x)dx.
$$
The operator $\nabla R_h:L^2_0(\Gamma_h)\rightarrow L^2(\Omega,\C^n)$ is defined on the whole $L^2_0(\Gamma_h)$ space, and the operator $F\mapsto AF$ is obviously bounded on $L^2(\Omega,\R^n)$. Moreover, if $\widetilde{u_h}\in W^{1,2}_0(\Gamma_h)$,
$$
\begin{array}{lll}
\displaystyle \int_{\Omega} A(x)\nabla\left(R_h\widetilde{u_h}\right)(x)\cdot \nabla\left(R_h\widetilde{u_h}\right)(x) dx & \geq & \delta \left\Vert \nabla \left(R_h\widetilde{u_h}\right)\right\Vert_{L^2(\Omega)}^2\\
& \gtrsim & \displaystyle \left\Vert R_h\widetilde{u_h}\right\Vert_{L^2(\Omega)}^2\\
& \gtrsim & \displaystyle \left\Vert \widetilde{u_h}\right\Vert_{L^{2}_0(\Gamma_h)}^2.
\end{array}
$$
As a consequence, using a construction due to Kato (\cite{kato}), there exists a unique maximal accretive operator $L_h$ on $L^2_0(\Gamma_h)$ such that, for all $\widetilde{u_h}\in {\mathcal D}(L_h)$ (here and after, ${\mathcal D}(L)$ stands for the domain of an operator $L$) and all $\widetilde{v_h}\in W^{1,2}_0(\Gamma_h)$,
$$
\langle L_h\widetilde{u_h},\widetilde{v_h}\rangle_{L^2_0(\Gamma_h)}=\int_{\Omega} A(x)\nabla\left(R_h\widetilde{u_h}\right)(x)\cdot \nabla\left(R_h\widetilde{v_h}\right)(x)dx.
$$
The following facts are easily checked:
\begin{pro} \label{galerkinlp}
For all $p\in (1,+\infty)$, the operator $L_h$ extends to a bounded operator from $W^{1,p}_0(\Gamma_h)$ to $W^{-1,p}(\Gamma_h)$. Moreover, $L_h$ is an isomorphism from $W^{1,2}_0(\Gamma_h)$ into $W^{-1,2}(\Gamma_h)$. The norms of $L_h$ and of its inverse (for $p=2$) do not depend on $h$. 
\end{pro}
{\bf Proof: } Let $p\in (1,+\infty)$ and $p^{\prime}\in (1,+\infty)$ such that $\frac 1p+\frac 1{p^{\prime}}=1$. Let $\widetilde{u_h}\in W^{1,p}_0(\Gamma_h)$. The map $\widetilde{v_h}\mapsto \int_{\Omega} A(x)\nabla\left(R_h\widetilde{u_h}\right)(x)\cdot \nabla\left(R_h\widetilde{v_h}\right)(x)dx$ is clearly a linear functional on $W^{1,p^{\prime}}_0(\Gamma_h)$. Moreover, for all $\widetilde{v_h}\in W^{1,p^{\prime}}_0(\Gamma_h)$,
$$
\begin{array}{lll}
\displaystyle \left\vert \langle L_h\widetilde{u_h},\widetilde{v_h}\rangle\right\vert & \lesssim & \displaystyle \left\Vert \nabla\left(R_h\widetilde{u_h}\right)\right\Vert_{L^p(\Omega)} \left\Vert \nabla\left(R_h\widetilde{v_h}\right)\right\Vert_{L^{p^{\prime}}(\Omega)}\\
& \lesssim & \displaystyle \left\Vert \nabla \widetilde{u_h}\right\Vert_{L^p(\Gamma_h)} \left\Vert \nabla\widetilde{v_h}\right\Vert_{L^{p^{\prime}}(\Gamma_h)}.
\end{array}
$$
This ends the proof of the first part of the claim. Assume now that $p=2$. For all bounded linear functional $T$ on $W^{1,2}_0(\Gamma_h)$, the Lax-Milgram theorem yields a unique function $\widetilde{u_h}\in W^{1,2}_0(\Gamma_h)$ such that, for all $\widetilde{v_h}\in W^{1,2}_0(\Gamma_h)$, $\langle L_h\widetilde{u_h},\widetilde{v_h}\rangle_{L^2_0(\Gamma_h)}=T(\widetilde{v_h})$, and one has $\left\Vert \widetilde{u_h}\right\Vert_{W^{1,2}_0(\Gamma_h)}\lesssim \left\Vert T\right\Vert_{W^{-1,2}(\Gamma_h)}$. This shows that $L_h$ is an isomorphism from $W^{1,2}_0(\Gamma_h)$ into $W^{-1,2}(\Gamma_h)$. Moreover, the norms of $L_h$ and $L_h^{-1}$ do not depend on $h$ because of Lemma \ref{geomgraph} and Proposition \ref{equivnorm}. \hfill\fin\par
\subsection{Conclusion of the proof}
Our next step is a perturbation argument, in order to show that $L_h$ is still an isomorphism from $W^{1,p}_0(\Gamma_h)$ into $W^{-1,p}(\Gamma_h)$ for $p$ close enough to $2$. Indeed, since $(P_{2,loc})$ holds, there exists $q_0<2$ such that, for $q_0<p_0<p<p_1<q'_0$, the $W^{1,p}(\Gamma)$ and $W^{-1,p}(\Gamma)$ form a complex interpolation scale (see Corollary \ref{complex1} and Theorem \ref{complex2}). We then combine this fact, Proposition \ref{galerkinlp} and the following general result (\cite{sneiberg}):
\begin{lem} \label{sneiberg} Let $X^s$, $Y^s$, $s\in [0,1]$ be two scales of complex interpolation Banach spaces. If $T:X^s \rightarrow Y^s$ is bounded for each $s\in [0,1]$, then the set of $s$ for which there exists
$C>0$ such that $\|Tf\|_{Y^s}\geq C\|f\|_{X^s}$ holds for all $f\in X^s$ is open.
\end{lem}
Thus, there exists $\varepsilon>0$ such that $L_h$ is an isomorphism from $W^{1,p}_0(\Gamma_h)$ into $W^{-1,p}(\Gamma_h)$ for all $p\in (2-\varepsilon,2+\varepsilon)$. \Bk Thus, we have established: 
\begin{pro} \label{lh}
The operator $L_h^{-1}$ maps  $W^{-1,p}(\Gamma_h)$ into $W^{1,p}_0(\Gamma_h)$ for $p\in (2,2+\varepsilon)$ with a norm independent of $h$.
\end{pro}

\medskip

\noindent Let us now conclude the proof of Theorem \ref{meyersgalerkin}. Take $f\in W^{-1,p}(\Omega)$ and define, for all $x\in \Gamma_h$,
$$
f_h(x):=
\left\{
\begin{array}{ll}
\displaystyle \frac 1{m(x)} \langle f,\varphi_x\rangle & \mbox{ when }x\in \Gamma_h\cap \Omega,\\
0 & \mbox{ when }x\in \Gamma_h\cap \partial\Omega.
\end{array}
\right.
$$
\Bk We claim that, if $u_h\in V^1_0({\mathcal T}_h)$ solves \eqref{galerkin} and $\widetilde{u_h}\in W^{1,2}_0(\Gamma_h)$ is defined by 
$$
\widetilde{u_h}(x)=
\left\{
\begin{array}{ll}
u_h(x) & \mbox{ when }x\in \Gamma_h\cap \Omega,\\
0 & \mbox{ when }x\in \Gamma_h\cap \partial\Omega,
\end{array}
\right.
$$
then 
\begin{equation} \label{lhuh}
L_h\widetilde{u_h}=-f_h.
\end{equation}
Indeed, $L_h\widetilde{u_h}\in W^{-1,2}(\Gamma_h)$, and \eqref{lhuh} exactly means that, for all $x\in \Gamma_h\cap \Omega$,
$$
-\langle L_h\widetilde{u_h},\widetilde{\varphi_x}\rangle= \langle f_h,\widetilde{\varphi_x}\rangle,
$$
where $\widetilde{\varphi_x}$ is the function in $W^{1,2}_0(\Gamma_h)$ defined by $\widetilde{\varphi_x}(y)=\delta_{xy}$ for all $y\in \Gamma_h$, so that $\varphi_x=R_h\widetilde{\varphi_x}$. But,
\begin{equation} \label{step1} 
\langle L_h\widetilde{u_h},\widetilde{\varphi_x}\rangle=m(x)L_h\widetilde{u_h}(x),
\end{equation}
and
\begin{equation} \label{step2}
\begin{array}{lll}
\displaystyle -\langle L_h\widetilde{u_h},\widetilde{\varphi_x}\rangle & = & \displaystyle - \int_{\Omega} A(y)\nabla\left(R_h\widetilde{u_h}\right)(y)\cdot \nabla\left(R_h\widetilde{\varphi_x}\right)(y)dy\\
& = & \displaystyle  -\int_{\Omega} A(y)\nabla u_h(y)\cdot \nabla \varphi_x(y)\\
& = &\Bk \displaystyle \langle f,\varphi_x\rangle \\ \Bk
& =& \displaystyle m(x)f_h(x)\\
& = & \displaystyle \langle f_h,\widetilde{\varphi_x}\rangle,
\end{array}
\end{equation}
where we used \eqref{galerkin} in the third equality and the definition of $f_h$ in the fourth one. Gathering \eqref{step1} and \eqref{step2}, we obtain \eqref{lhuh}. It therefore follows from \eqref{lhuh} and Proposition \ref{lh} that, for all $p\in (2,2+\varepsilon)$,
\begin{equation} \label{estimtildeuh}
\left\Vert \widetilde{u_h}\right\Vert_{W^{1,p}(\Gamma_h)} \lesssim \left\Vert f_h\right\Vert_{W^{-1,p}(\Gamma_h)}.
\end{equation}
Proposition \ref{equivnorm} shows that, \Bk for all $p\in (2,2+\varepsilon)$,
\begin{equation} \label{equivuhtildeuh}
\Bk \left\Vert u_h\right\Vert_{W^{1,p}(\Omega)}\sim \left\Vert \widetilde{u_h}\right\Vert_{W^{1,p}(\Gamma_h)}. \Bk
\end{equation}
Moreover,
\begin{equation} \label{ftildefh}
\Bk \left\Vert f_h\right\Vert_{W^{-1,p}(\Gamma_h)}\leq C\left\Vert f\right\Vert_{W^{-1,p}(\Omega)}.
\end{equation}
Indeed, let $v$ be a function in $\Gamma_h$ with finite support. Then,
$$
\begin{array}{lll}
\displaystyle \left\vert \langle f_h,v\rangle\right\vert & = &\displaystyle \left\vert  \sum_{x\in \Gamma_h} \langle f,\varphi_x\rangle v(x)\right\vert \\
& = & \displaystyle \left\vert  \langle f,\sum_{x\in \Gamma_h} v(x)\varphi_x\rangle \right\vert \\
& \leq & \displaystyle \left\Vert f\right\Vert_{W^{-1,p}(\Omega)} \left\Vert \sum_{x\in \Gamma_h} v(x)\varphi_x \right\Vert_{W^{1,p^{\prime}}_0(\Omega)} \\
&=  & \displaystyle  \left\Vert f\right\Vert_{W^{-1,p}(\Omega)} \left\Vert R_hv\right\Vert_{W^{1,p^{\prime}}_0(\Omega)} \\
& \lesssim & \displaystyle  \left\Vert f\right\Vert_{W^{-1,p}(\Omega)} \left\Vert v\right\Vert_{W^{1,p^{\prime}}_0(\Gamma_h)},
\end{array}
$$
which proves \eqref{ftildefh}. 
\Bk Gathering \eqref{estimtildeuh}, \eqref{equivuhtildeuh} and \eqref{ftildefh} yields \eqref{estimuh}.\par
\noindent \Bk For the proof of \eqref{convuh}, it is enough to interpolate between \eqref{estimuh} and the fact that $\left\Vert u-u_h\right\Vert_{W^{1,2}(\Omega)}\rightarrow 0$ when $h\rightarrow 0$. \hfill\fin\par
\noindent Let us now prove Proposition \ref{optimal}. Let $\varepsilon>0$. By \cite{meyers}, there exists a uniformly elliptic matrix $A\in L^{\infty}(\Omega)$ with the following property: for all $p>2+\varepsilon$, there exists $f\in W^{-1,p}(\Omega)$ such that the solution of $-\mbox{div}(A\nabla u)=f$ in $W^{1,2}_0(\Omega)$ does not belong to $W^{1,p}(\Omega)$. Assume now that the $u_h$ are uniformly bounded in $W^{1,p}(\Omega)$. Up to a subsequence, $\nabla u_h$ converges weakly in $L^p(\Omega)$ to some vector field $V$. But since $\nabla u_h\rightarrow \nabla u$ in the $L^2$ norm, it follows that $\nabla u=V$ and that $\nabla u\in L^p(\Omega)$, a contradiction. \hfill\fin\par
\noindent As explained in the introduction, Corollary \ref{holdgalerkin} follows at once from Theorem \ref{meyersgalerkin} by Sobolev embeddings. Finally, for Corollary \ref{improvedgalerkin}, it is well-known (see \cite{bs}, Chapter 5) that, in this situation, 
\begin{equation} \label{speed}
\left\Vert u-u_h\right\Vert_{W^{1,2}(\Omega)}\leq Ch\left\Vert f\right\Vert_{L^2(\Omega)}.
\end{equation}
By the results of \cite{galmon} stated in the introduction, and since $L^2(\Omega)\hookrightarrow W^{-1,2+\varepsilon}(\Omega)$, one has 
\begin{equation} \label{sobnormu}
\left\Vert u\right\Vert_{W^{1,2+\varepsilon}(\Omega)}\leq C\left\Vert f\right\Vert_{L^2(\Omega)}.
\end{equation}
so that
\begin{equation} \label{sobu-uhp}
\left\Vert u-u_h\right\Vert_{W^{1,2+\varepsilon}(\Omega)}\leq C\left\Vert f\right\Vert_{L^2(\Omega)}.
\end{equation}
Interpolating between \eqref{sobu-uhp} and \eqref{speed} gives \eqref{convuhsobolev} with
$$
\theta:=\frac{\frac 1p-\frac 1{2+\varepsilon}}{\frac 12-\frac 1{2+\varepsilon}}.\hfill\fin
$$
\Bk \section{Estimates for general second elliptic operators on graphs}\label{ellip}
In the present section, we prove estimates in $L^{\infty}$ and H\"older spaces for general second order elliptic on graphs, assuming suitable geometric properties on the graph $\Gamma$, which is assumed to be unbounded.
\subsection{Definition of the operators}
Assume that, for all $x\sim y$, a coefficient $c_{xy}\in \C$ is given, and that there exist $C_{\infty}>0$ such that
\begin{equation} \label{bounded}
\left\vert c_{xy}\right\vert \leq C_{\infty}\ \forall x\sim y
\end{equation}
and $\delta>0$ such that, for all $u\in W^{1,2}(\Gamma)$, 
\begin{equation} \label{elliptic}
\mbox{Re }\sum_{x\sim y} c_{xy}\left\vert du(x,y)\right\vert^2\mu_{xy}\geq \delta \sum_{x\sim y} \left\vert du(x,y)\right\vert^2\mu_{xy}.
\end{equation}
We associate to these coefficients an operator, which is the discrete version of second order uniformly elliptic operators in divergence form in $\R^n$. To that purpose, we use again a classical construction of maximal accretive operators due to Kato (\cite{kato}). The operator $d:L^2(\Gamma)\rightarrow L^2(E)$ is densely defined (see Remark \ref{denselp}) and closed, and the operator $F\mapsto cF$ is clearly $L^2(E)$ bounded by \eqref{bounded}. Since \eqref{elliptic} holds, there exists a unique maximal accretive operator $L$ on $L^2(\Gamma)$ such that, for all $u\in {\mathcal D}(L)$ and all $v\in W^{1,2}(\Gamma)$,
\begin{equation} \label{defL}
\sum_{x\in \Gamma} Lu(x)\overline{v(x)}m(x)=\sum_{(x,y)\in E} c_{xy} du(x,y)\overline{dv(x,y)}\mu_{xy}
\end{equation}
(recall that ${\mathcal D}(L)$ stands for the domain of $L$).
\begin{rem}
Note that assumption \eqref{elliptic} is satisfied in particular when the coefficients $c_{xy}$ are real-valued and satisfy $c_{xy}\geq \delta$ whenever $x\neq y$. However, we do {\bf not} make this assumption in the sequel, and all that is needed is \eqref{elliptic}. 
\end{rem}
We say that $L$ is a second order uniformly elliptic operator on $L^2(\Gamma)$ with ellipticity constants $C_{\infty}$ and $\delta$. Recall that ${\mathcal D}(L)$ is dense in $L^2(\Gamma)$ and that $L$ generates a holomorphic semigroup on $L^2(\Gamma)$. Note finally that, for all $u\in {\mathcal D}(L)$,
\begin{equation} \label{ellipticL}
\mbox{Re }\langle Lu,u\rangle \geq c\delta \left\Vert \nabla u\right\Vert_2^2.
\end{equation} 
\begin{rem} \label{Lsobolev}
Let $1<p<+\infty$. Assume that $u\in {\mathcal D}(L)\cap W^{1,p}(\Gamma)$ and $v\in W^{1,2}(\Gamma)\cap W^{1,p^{\prime}}(\Gamma)$. From the definition of $L$, we deduce that
$$|\langle Lu,v \rangle| \leq C \|\nabla u\|_{L^p(\Gamma)} \|\nabla v\|_{L^{p'}(\Gamma)} .
$$
This and Remark \ref{densesobol} show that $Lu$ extends to a bounded  linear form on $W^{1,p'}(\Gamma)$, with norm less or equal to $C \| u\|_{W^{1,p}(\Gamma)}$ . Thus $L$ is a bounded operator from $W^{1,p}(\Gamma)$ to  $W^{-1,p}(\Gamma)$, with bound less or equal to $C$.
\end{rem}

\subsection{The estimates on $L$}
Our first estimates deal with the resolvent of $L$. Let $\omega\in \left(0,\frac{\pi}2\right)$ be such that $L$ is $\omega$-accretive on $L^2(\Gamma)$, $\mu \in \left(\frac{\pi}{2},\pi-\omega\right)$ and
$$
\Sigma_{\mu}:=\left\{z\in \C;\ z=0\mbox{ or }\left\vert \mbox{arg }z\right\vert<\mu\right\}.
$$
Let $f\in L^2(\Gamma)$ and $\lambda\in \Sigma_{\mu}$. We are interested in weak solutions of the equation
\begin{equation} \label{weaksol}
Lu+\lambda u=f.
\end{equation}
By ``weak solution'', we mean a function $u\in W^{1,2}(\Gamma)$ such that, for all $v\in W^{1,2}(\Gamma)$,
$$
\sum_{(x,y)\in E} c_{xy}du(x,y)\overline{dv(x,y)}\mu_{xy}+\lambda\sum_{x\in \Gamma} u(x)\overline{v(x)}m(x)=\sum_{x\in \Gamma} f(x)\overline{v(x)}m(x).
$$
We establish that, if  $(D)$, $(L_2)$ and $(P_2$) hold, given $f\in L^2(\Gamma)$ and $\lambda\in \Sigma_{\mu}$, \eqref{weaksol} has a unique weak solution in $W^{1,2}(\Gamma)$, and that this solution actually belongs to $C^{\eta}(\Gamma)$ for some $\eta>0$:
\begin{theo} \label{resolventholder}
Let $\Gamma$ be a graph satisfying $(D)$, $(L_2)$ and $(P_2)$. Let $L$ be a second order uniformly elliptic operator on $L^2(\Gamma)$  with ellipticity constants $C_{\infty}$ and $\delta$. Let $\omega$ be such that $L$ is $\omega$-accretive on $L^2(\Gamma)$ and $\mu \in \left(\frac{\pi}{2},\pi-\omega\right)$. Then, for all $f\in L^2(\Gamma)$ and all $\lambda\in \Sigma_{\mu}$, there exists a unique solution $u\in W^{1,2}(\Gamma)$ of \eqref{weaksol}. Moreover, there exist $C,\eta>0$ only depending on the constants in  $(D)$, $(L_2)$ and $(P_2)$ and the constants of ellipticity of $L$ such that
$$
\left\Vert u\right\Vert_{L^{\infty}(\Gamma)}\leq C\left\vert \lambda\right\vert^{-\frac 12} \left\Vert f\right\Vert_{L^2(\Gamma)}
$$
and
$$
\left\vert u\right\vert_{\dot{C}^{\eta}(\Gamma)}\leq C\left\vert \lambda\right\vert^{\frac{\eta}2} \left\vert \lambda\right\vert^{-\frac 12} \left\Vert f\right\Vert_{L^2(\Gamma)}.
$$
\end{theo}
Theorem \ref{resolventholder} is a version of the De Giorgi regularity theorem for elliptic equations in this context.\par
\noindent We also obtain estimates on the kernel of the semigroup generated by $L$. Before  stating the Theorem, we introduce  the following notations: for all $x,\, y \in \Gamma$, $h^*_{x\rightarrow y}$ is equal to $0$ if $x=y$ and
is the supremum of the weights of edges where at least one of the two vertices  belongs to the ball $B(x,d(x,y))$, when $x \neq y$. We then define
$$h^*_{xy}= \min(h^*_{x\rightarrow y}, h^*_{y\rightarrow x}).$$

\begin{theo} \label{estimationsholderiennes} Let $\Gamma$ be a graph satisfying $(D)$, $(L_2)$ and $(P_2)$. Let $L$ be a second order uniformly elliptic operator on $L^2(\Gamma)$  with ellipticity constants $C_{\infty}$ and $\delta$. We then have:

\begin{itemize}
\item[$1.$] 
The operator $L$ generates a holomorphic semi-group $(e^{-tL})_{t>0}$ on $L^2(\Gamma)$ which has a kernel, denoted by $K_{t}(x,y)$ in the sequel. For all $u \in L^{2}(\Gamma)$ and all $x \in \Gamma$,
$$
(e^{-tL}u)(x) =\sum_{y\in \Gamma}K_t(x,y)u(y).
$$
Moreover, there exist $C, \,C',\beta>0$ depending only of the geometric constants and that of ellipticity, such that for every $(x,y)\in \Gamma^2$ and $t>0$ we have
\begin{itemize}
\item [a.]if $t\leq C'h^*_{xy}d(x,y) $,
$$|K_{t}(x,y)|\leq \frac{C}{t} e^{-\beta \frac{d(x,y)}{h^*_{xy}}};
$$
\item[b.] if $t\geq C'h^*_{xy}d(x,y) $,
$$|K_{t}(x,y)|\leq \frac{C}{t} e^{-\beta \frac{d^2(x,y)}{t}}.
$$
\end{itemize}
\item[$2.$] 
There exist constants $C''$ and $\eta>0$ depending only on the geometric constants and that of ellipticity such that for every $t>0$ and for every $x,x^{\prime},y\in \Gamma$, we have
$$
|K_t(x,y)-K_t(x^{\prime},y)|\leq \frac{C''}{t}\left(\frac{d(x,x^{\prime})}{\sqrt{t}}\right)^{\eta}.
$$
\end{itemize}
\end{theo}
\subsection{Proofs}
We first focus on Theorem \ref{resolventholder}. We begin with the following observation:
\begin{pro} Let $\Gamma$ as above and $L$ a second order elliptic operator on $\Gamma$, with ellipticity constants $C_{\infty}$ and $\delta$.
Then, for all $\lambda\in \C$, $L+\lambda I$ is a continuous operator from $W^{1,2}(\Gamma)$ to $W^{-1,2}(\Gamma)$. Moreover, when $\lambda=1$, $L+I$ is an isomorphism from $W^{1,2}(\Gamma)$ onto $W^{-1,2}(\Gamma)$ and the norm of its inverse is bounded by a constant only depending on $\delta$.
\end{pro}
\textbf{Proof: } That $L+\lambda I$ is bounded from $W^{1,2}(\Gamma)$ to $W^{-1,2}(\Gamma)$ was already seen in Remark \ref{Lsobolev} above. The Lax-Milgram theorem, applied with the sesquilinear form
$$
{\mathcal B}(u,v):=\sum_{(x,y)\in E}c_{x,y}du(x,y)\overline{dv(x,y)}\mu_{xy}+\sum_{x\in \Gamma} u(x)\overline{v(x)}m(x)
$$
which is clearly continuous and coercive on $W^{1,2}(\Gamma)$ thanks to \eqref{ellipticL}, yields the invertibility of $L+I$ and the bound on the norm of the inverse. \hfill\fin\par

\medskip

\noindent Relying on Lemma \ref{sneiberg} again and arguing as in Proposition \ref{lh}, we prove:
\begin{pro}\label{appsneiberg} Assume that $\Gamma$ satisfies (\ref{D}) and $(P_{2})$ and $L$ is an elliptic operator as above. Then there exists $\varepsilon>0$ only depending on the geometric constants of $\Gamma$ and the ellipticity constants of $L$ such that $L+I$ is invertible from $W^{1,p}$ to $W^{-1,p}$ for $p\in (2-\varepsilon,2+\varepsilon)$.
\end{pro}

\subsection{Proof of the estimates on $L$} \label{proofestim}
We prove here Theorems \ref{resolventholder} and \ref{estimationsholderiennes}.
\subsection{Proof of Theorem \ref{resolventholder}}
Let $f\in L^2(\Gamma)$. Since $L+I$ is an isomorphism from $W^{1,2}(\Gamma)$ onto $W^{-1,2}(\Gamma)$, there exists a unique $u\in W^{1,2}(\Gamma)$ such that \eqref{weaksol} is satisfied.  Moreover, we have seen in Proposition \ref{appsneiberg} that, for some $\varepsilon>0$ $L+I$ is an isomorphism from $W^{1,p}(\Gamma)$ onto $W^{-1,p}(\Gamma)$ for $p\in (2-\varepsilon,2+\varepsilon)$. Let $p\in (2,2+\varepsilon)$. Propositions \ref{sobolev}  and \ref{holder} applied with $q=\sigma=2$ yield
$$
L^2(\Gamma) \hookrightarrow W^{-1,p}(\Gamma)
$$
and
$$
W^{1,p}(\Gamma)\hookrightarrow C^{\eta}(\Gamma).
$$
Therefore $(L+I)^{-1}$ maps  $L^2(\Gamma)$ into $C^{\eta}(\Gamma)$ with $\eta=1-\frac{2}{p}$. \par
\noindent With an analogous argument, $L+I$ can be replaced by $L+\lambda I$ for $\lambda\in \Sigma_{\mu}$ with $|\lambda|=1$:
\begin{equation} \label{R1}
\|(L+\lambda)^{-1}u\|_{\infty}\leq C \|u\|_{L^2(\Gamma)}
\end{equation}
and
\begin{equation}\label{R2}
\left\vert (L+\lambda)^{-1}u\right\vert_{\dot{C}^{\eta}}\leq C \|u\|_{L^2(\Gamma)}.
\end{equation}
We now claim that there is a constant $C>0$ such that, for all $\lambda\in \Sigma_{\mu}\setminus \left\{0\right\}$,
$$
\|(L+\lambda)^{-1}u\|_{\infty}\leq C |\lambda|^{-\frac{1}{2}}\|u\|_{L^2(\Gamma)}
$$
and
$$
|(L+\lambda)^{-1}u|_{\dot{C}^{\eta}}\leq C  |\lambda|^{\frac{\eta}{2}-\frac{1}{2}}\|u\|_{L^2(\Gamma)}.
$$
We treat the case where $\lambda >0$ using a scaling argument. A similar proof can be done for a general $\lambda \in \Sigma_{\mu}$.  \par
\noindent The scaling argument is as follows. For all $\alpha>0$, let $\Gamma_{\alpha}$ be the set $\Gamma$, endowed with the weights $\alpha h_{xy}$ and the measures $\alpha^2\mu_{xy}$. Note that, if $\nabla_{\alpha}$ stands for the gradient in $\Gamma_{\alpha}$, one has $\nabla_{\alpha}f=\frac 1{\alpha}\nabla f$. As a consequence, it is easily checked that, endowed with the distance $d_{\alpha}$ defined by the weights $\alpha h_{xy}$ and the measure $\alpha^2m$, the graph $\Gamma_{\alpha}$ satisfies the assumptions $(D)$, $(P_2)$ and $(L_2)$ with the same constants as $\Gamma$. \par
\noindent Consider now the operator $L_{\alpha}$ on $\Gamma_{\alpha}$ given by the coefficients $c_{xy}$. A straightforward computation yields that $L_{\alpha}$ has the same ellipticity constants on $\Gamma_{\alpha}$ as $L$ on $\Gamma$. \par
\noindent Let $f\in L^2(\Gamma)$ and $\lambda>0$ be given. Define $g:=\frac 1{\lambda}f$. Applying the conclusion of Theorem \ref{resolventholder} with $\Gamma_{\sqrt{\lambda}}$ and the operator $L_{\sqrt{\lambda}}$, one obtains that there exists a unique function $u\in W^{1,2}(\Gamma_{\sqrt{\lambda}})$ such that
$$
(L_{\sqrt{\lambda}}+I)u=g
$$
on $\Gamma_{\sqrt{\lambda}}$, and 
$$
\left\Vert u\right\Vert_{C^{\eta}(\Gamma_{\sqrt{\lambda}})}\leq C\left\Vert g\right\Vert_{L^2(\Gamma_{\sqrt{\lambda}})}=\frac C{\sqrt{\lambda}}\left\Vert f\right\Vert_{L^2(\Gamma)}.
$$
Therefore, $(L+\lambda I)u=f$ on $\Gamma$. Moreover,
$$
\left\Vert u\right\Vert_{L^{\infty}(\Gamma)}\leq \frac C{\sqrt{\lambda}}\left\Vert f\right\Vert_{L^2(\Gamma)},
$$
and
$$
\left\vert u\right\vert_{\dot{C}^{\eta}(\Gamma)}\leq C\lambda^{\frac{\eta-1}2} \left\Vert f\right\Vert_{L^2(\Gamma)},
$$
which ends the proof of Theorem \ref{resolventholder}. \hfill\fin
\subsection{Proof of Theorem \ref{estimationsholderiennes}}
The proof of Theorem \ref{estimationsholderiennes}
is analogous to that of Proposition 28, p.59 in \cite{asterisque} in the Euclidean case, and we only indicate the main differences.\par
\noindent Let us first recall that, for all $t>0$, 
\begin{equation}\label{cauchy}
e^{-tL}=\frac{1}{2\pi i}\int_{\gamma}e^{z\lambda}(L+\lambda)^{-1}d\lambda
\end{equation}
where $\gamma$ is made of two rays $\gamma_{\pm}:=\left\{re^{\pm i\theta};\ r\geq \frac 1t\right\}$ and of the arc $\gamma_0:=\left\{\frac 1t e^{i\sigma};\ \vert \sigma\vert\leq \theta\right\}$ where $\theta\in ]\frac{\pi}2,\pi[$, and $\gamma$ is described counterclockwise.\par
\noindent Combining \eqref{cauchy} with Theorem \ref{resolventholder}, we obtain
$$
\|e^{-tL}u\|_{L^{\infty}(\Gamma)}\leq Ct^{-\frac{1}{2}}\|u\|_{L^2(\Gamma)},\qquad t>0
$$
and
$$
|e^{-tL}u|_{\dot{C}^{\eta}}\leq Ct^{-(\frac{1}{2}+\frac{\eta}{2})}\|u\|_{2} ,\qquad t>0.
$$
We finish the proof as in \cite{rey}, Chapter 5, section 5.2, see also \cite{asterisque}, Chapter 1. \hfill\fin\par

\end{document}